\documentclass[12pt,a4paper]{amsart}

\usepackage{amsmath}
\usepackage{amssymb}
\usepackage{amscd}

\newtheorem{thm}{Theorem}[section]

\newtheorem{lem}[thm]{Lemma}

\newtheorem{cor}[thm]{Corollary}
\theoremstyle{remark}
\newtheorem{example}[thm]{Example}
\newtheorem{examples}[thm]{Examples}
\newtheorem{remark}[thm]{Remark}

\def\N{{\mathbb N}}
\def\A{{\mathbb A}}
\def\Q{{\mathbb Q}}
\def\T{{\mathbb T}}
\def\H{{\mathcal{H}}}
\def\Z{{\mathbb Z}}

\newcommand{\End}{\operatorname{End}}
\newcommand{\Ind}{\operatorname{Ind}}
\newcommand{\Aut}{\operatorname{Aut}}
\def\dashind{\operatorname{\!-Ind}}

\numberwithin{equation}{section}

\begin{document}

\title[Representations of Hecke algebras]{Representations of Hecke
    algebras\\ and dilations of semigroup crossed products}
\author[Nadia S. Larsen]{Nadia S. Larsen}
\address{Department of Mathematics, University of Copenhagen,
Universitetsparken 5, DK-2100 Copenhagen \O, Denmark.}
\email{nadia@math.ku.dk}
\thanks{The first author and her research were supported by the Danish
  Natural Science Research Council, and the research was also supported by
the Australian Research Council.}
\author[Iain Raeburn]{Iain Raeburn}
\address{Department of Mathematics, University of Newcastle, NSW 2308,
Australia.}
\email{iain@maths.newcastle.edu.au}
\date {8 June 2001}

\begin{abstract}
We consider a family of Hecke $C^*$-algebras which can be realised as crossed
products by semigroups of endomorphisms. We show by dilating representations of
the semigroup crossed product that the category of representations of the Hecke
algebra is equivalent to the category of continuous unitary representations of
a totally disconnected locally compact group.
\end{abstract}

\maketitle

Suppose that $M$ is a subgroup of a group $\Gamma$ such that
$|M\cap\gamma^{-1}M\gamma|$ is finite for every $\gamma\in
\Gamma$; we say that $M$ is an almost normal subgroup or that
$(\Gamma,M)$ is a Hecke pair. The Hecke algebra
$\H(\Gamma,M)$ is a convolution
$*$-algebra of functions on the double coset space $M\backslash\Gamma/M$,
which can be represented in the commutant of the quasi-regular
representation $\Ind_M^\Gamma 1$ \cite{K}; $C^*$-algebraic completions of
Hecke algebras played a fundamental role in the analysis by Bost and
Connes of phase transitions in number theory
\cite{BC}. Several authors have since investigated classes of $C^*$-Hecke
algebras which can be realised as semigroup crossed products \cite{LR3,
ALR, B, LL}, and these realisations have provided valuable insight to
the work of Bost and Connes \cite{L2, N}.

A recent theorem of Hall \cite{H} asserts
that, for a large class of Hecke pairs $(\Gamma, M)$, the category of
nondegenerate representations of the Hecke algebra $\H(\Gamma, M)$ is
equivalent to the category $\mathcal{C}(\Gamma,M)$ of unitary
representations of $\Gamma$ which are generated by their $M$-fixed
vectors. From an operator-algebraic point of view, $\mathcal{C}(\Gamma,M)$
seems an unusual category: it is not obvious, for example, whether
it is the category of representations of some familiar $C^*$-algebra.
Here we consider a class of Hecke algebras $\H(\Gamma, M)$ which can be
realised as semigroup crossed products, and show that for these pairs
$(\Gamma,M)$, a representation of $\Gamma$ is in $\mathcal{C}(\Gamma,M)$
precisely when it has been dilated from a representation of the
corresponding semigroup crossed product. We then use the
dilation-extension theory of Laca \cite{L} to identify a single locally
compact group $\Gamma_\infty$ whose category of continuous unitary
representations is equivalent to $\mathcal{C}(\Gamma,M)$.

We begin in \S\ref{Hallsthm} by 
describing a family $(\Gamma,M)$ of Hecke pairs in which $\Gamma$
is a semi-direct product $N\rtimes G$, $M$ is a normal subgroup of $N$,
and $G$ has the form $S^{-1}S$ for some subsemigroup $S$. It is shown in
\cite{LL}, generalising results in \cite{B} and \cite{ALR}, that the associated
Hecke algebra
$\H(\Gamma, M)$ is isomorphic to a semigroup crossed product
$C^*(N/M)\rtimes_\alpha S$. In Theorem~\ref{Hallviadilation}, we show that for
each covariant representation $(\pi,V)$ of $(C^*(N/M),S,\alpha)$ and each
minimal unitary dilation $U$ of $V$, there is a representation of $N\rtimes G$
on
$\H_U$ which is generated by its $M$-fixed vectors; more importantly, we
show that they all arise this way.

Laca's analysis of dilations now suggests that we should study the minimal
automorphic dilation of the semigroup system $(C^*(N/M),S,\alpha)$, which is an
ordinary dynamical system $(B,G,\beta)$ such that $C^*(N/M)\rtimes_\alpha
S$ embeds naturally as a full corner in $B\rtimes_\beta G$. We take
for $B$ the $C^*$-algebra of a totally disconnected locally compact
completion
$N_\infty$ of
$N$, so that $B\rtimes_\beta G$ is the group $C^*$-algebra of a
semi-direct product $\Gamma_\infty=N_\infty\rtimes G$ which is a
completion of
$N\rtimes G$. For the Hecke algebra of Bost and Connes, our construction
yields as $N_\infty$ the group $\A_f$ of finite adeles, as in
\cite[\S3]{L}. 

Our main Theorem~\ref{idcategory} says that $\mathcal{C}(\Gamma,M)$ is
equivalent to the category of continuous unitary representations of
$N_\infty\rtimes G$, and that this equivalence is implemented by Laca's
dilation-extension process. To prove this last point, we need a
categorical version of \cite[Theorem~2.4]{L}, which we provide in the
Appendix. The main idea is that the dilation-extension process can be
viewed as the induction process associated to a Morita equivalence, and
hence is automatically an equivalence of categories.

While we have been working on this project, there have been several other
developments in the area. Tzanev \cite[Theorem~1.17]{Tz}, and independently
Gl{\" o}ckner and Willis \cite{GW}, have  used a construction of Schlichting
\cite{S} to associate to each Hecke pair $(\Gamma, M)$ a locally compact
completion $(\overline{\Gamma},\overline{M})$ which yields the
same Hecke algebra. The completion $\overline{\Gamma}$
is totally disconnected, and $\overline{M}$ is a compact open subgroup; a
uniqueness theorem for such completions
\cite[Proposition~1.16]{Tz} implies that, in our setting, $\overline{\Gamma}$ is
isomorphic to
$N_\infty\rtimes G$. We suspect, however, that our construction of $N_\infty$
as an inverse limit  may be helpful in applications (see \S\ref{exs}). Our
interest in this subject was aroused by conversations with Steve Kaliszewski,
Magnus Landstad and John Quigg, and they have now obtained Hall's equivalence
as a Morita equivalence associated to a corner in
$C^*(\overline{\Gamma})$ \cite{KLQ}. It is interesting that the Ore condition
$G=S^{-1}S$, which arises here because it is a crucial hypothesis in both
\cite{LL} and
\cite{L}, also appears in the analyses of \cite{H} and
\cite{KLQ}.

\subsection*{Acknowledgements} 

We thank Steve Kaliszewski, Magnus Landstad and
John Quigg for many helpful conversations.

\subsection*{Conventions}

We view the $C^*$-algebra of a group $G$ as the $C^*$-algebra
$C^*(G)$ generated by a universal unitary representation $g\mapsto
\delta_g$: when $G$ is discrete, the $\delta_g$ are unitary elements of
$C^*(G)$ which span a dense subspace of $C^*(G)$, and when
$G$ is locally compact, the $\delta_g$ are unitary multipliers of
$C^*(G)$,
$\delta:G\to UM(C^*(G))$ is strictly continuous, and the elements
$f=\int_G f(g)\delta_g\,dg$ for $f\in C_c(G)$ span a dense subspace of
$C^*(G)$. The map
$\pi\mapsto
\pi\circ \delta$ sets up a one-to-one correspondence between the
nondegenerate representations of $C^*(G)$ and the continuous unitary
representations of $G$; we write $\pi_W$ for the representation of $C^*(G)$
associated to $W:G\to U(\H)$, so that $W=\pi_W\circ \delta$.

If $\psi:G\to \Aut N$ is a homomorphism, we realise the semi-direct
product
$N\rtimes_\psi G$ as the set $N\times G$ with
$(m,g)(n,h):=(m\psi_g(n),gh)$. The action $\psi$ of $G$ on $N$ induces an
action $\psi_*$ of $G$ on $C^*(N)$ such that
$(\psi_*)_g(\delta_n)=\delta_{\psi_g(n)}$. If $W:N\to U(\H)$ and $U:G\to
U(\H)$ are unitary representations, then $(\pi_W,U)$ is a covariant
representation of $(C^*(N),G,\psi_*)$ if and only if
$U_gW_nU_g^*=W_{\psi_g(n)}$, and then $W\rtimes U:(n,g)\mapsto W_nU_g$
is a unitary representation of $N\rtimes_\psi G$. All
representations of
$N\rtimes_\psi G$ arise this way: indeed, there is a natural isomorphism
of $C^*(N\rtimes_\psi G)$ onto the crossed product
$C^*(N)\rtimes_{\psi_*}G$ which carries $\pi_{W\rtimes\, U}$ into the integrated
form $\pi_W\times U$ of $(\pi_W,U)$. All this works for discrete groups
and for continuous actions on locally compact groups.

\section{Hall's induction process as a dilation}\label{Hallsthm}

We begin by establishing our notation. Throughout $S$ will be an Ore
semigroup; this means that $S$ is a cancellative semigroup with  $Ss\cap
St\not=\emptyset$ for all $s,t\in S$, and implies that $S$ is a subsemigroup
of an enveloping group $G$ such that $G=S^{-1}S:=\{s^{-1}t:s,t\in S\}$.  We
suppose that we have an action
$\psi$ of $G$ by automorphisms of another group $N$, and that
there is a normal subgroup $M$ of $N$ such that $M\subset
\psi_s(M)$ and $|M/\psi_s^{-1}(M)|<\infty$ for $s\in S$, and such that
\begin{equation}
\textstyle{\bigcap_{s \in S}\psi_s^{-1}(M)=\{e\}.} \label{j_injective}
\end{equation}
It is shown in \cite{LL} that $M$ is an almost normal subgroup of
the semi-direct product $\Gamma:=N\rtimes_\psi G$, and thus gives rise to
a Hecke algebra $\H(\Gamma,M)$. The example to bear in mind is
$S=\N^*$, $G=\Q_+^*$, $N=\Q$, $\psi_g(r)=g^{-1}r$, and $M=\Z$, which yields the
Hecke algebra $\H(\Q\rtimes\Q_+^*,\Z)$ of Bost and Connes \cite{BC}.

It was also proved in \cite{LL} that there is an action $\alpha$ of $S$
by injective endomorphisms of $C^*(N/M)$ such that
\[
\alpha_s(\delta_{nM})= \frac 1{|M:\psi_s^{-1}(M)|}\sum_{\{mM\in
N/M\,:\,
\psi_s^{-1}(m)M=nM\}}\delta_{mM}.
\]
Every such action $\alpha:S\to \End A$ has a
crossed product $A\rtimes_\alpha S$, which is universal for
covariant representations $(\pi,V):(A,S)\to B(\H)$ consisting of a
representation $V$ of $S$ by isometries on $\H$ and a nondegenerate
representation $\pi$ of $A$ on $\H$ satisfying
$\pi(\alpha_s(a))=V_s\pi(a)V_s^*$ for
$a\in A$, $s\in S$. These \emph{semigroup crossed products} are
studied in \cite{LR2} and
\cite{M}, for example, and that every system of our kind has one is proved in
\cite{L}. A theorem in
\cite{LL} says that the Hecke algebra
$\H(\Gamma,M)$ is isomorphic to the semigroup crossed product
$C^*(N/M)\rtimes_\alpha S$; in our motivating example, we recover the
isomorphism $\H(\Q\rtimes\Q_+^*,\Z)\cong C^*(\Q/\Z)\rtimes_\alpha \N^*$
used in \cite{LR3},
\cite{L2}, \cite{LR4} and \cite{N}.

We can study isometric representations of a semigroup by dilating them
to unitary representations of an enveloping group.
For abelian semigroups, this construction has been around for many
years, and Laca has recently extended it to Ore semigroups~\cite{L}. He
proved that for every isometric representation $V:S\to B(H)$, there are a
Hilbert space $\H$ containing $H$ as a closed subspace and a unitary
representation
$U:G:=S^{-1}S\to B(\H)$ such that $V_s=U_s|_H$  and
$\bigcup_{s\in S}U_s^*H$ is dense in $\H$ \cite[Theorem~1.4]{L}; up to
unitary equivalence, there is exactly one such
\emph{minimal unitary dilation} $U$.

Our first theorem shows how
dilating representations of the semigroup system $(C^*(N/M),S,\alpha)$
yields Hall's category of representations of $N\rtimes_\psi G$ generated
by their
$M$-fixed vectors.

\begin{thm}\label{Hallviadilation}
Let $(\pi, V)$ be a covariant representation of $(C^*(N/M),S,\alpha)$ on
$H$, and suppose that $U:G\to U(\H)$ is a minimal unitary dilation of
$V$. Then there is a unitary representation $W$ of $N$ on $\H$ such that
\begin{equation}\label{defW}
W_nU_s^*h=U_s^*\pi(\delta_{\psi_s(n)M})h\ \text{ for $h\in H$ and $s\in
S$};
\end{equation}
in particular, $H\subset \H^M:=\{h\in \H:W_mh=h\text{ for all $m\in
M$}\}$. The pair $(W,U)$ defines a representation
$W\rtimes U$ of
$N\rtimes_\psi G$ which is generated by its space $\H^M$ of $M$-fixed vectors,
and we can recover $\pi$ as  $\pi_Y$, where $Y_{nM}:=W_n|_H$.

Conversely, suppose $X\rtimes U:N\rtimes_\psi G\to U(\H)$ is a
unitary representation which is generated by its space $\H^M$ of $M$-fixed
vectors. Then
$V:s\mapsto U_s|_{\H^M}$ is an isometric representation of $S$, and $U$
is a minimal unitary dilation of $V$. With $Y_{nM}:=X_n|_{\H^M}$, the
pair $(\pi_Y,V)$ is a covariant representation of $(C^*(N/M),S,\alpha)$
on $\H^M$, and applying the construction of the previous paragraph to
$(\pi_Y,V)$ gives $W=X$.
\end{thm}

The proof of this theorem will occupy the rest of the section. For the
first part, suppose $(\pi,V)$ is a covariant representation of
$(C^*(N/M),S,\alpha)$ on $H$. We write $\pi=\pi_Y$, and the covariance
relation translates to
\begin{equation}
V_sY_{nM}V_s^*=\frac{1}{\vert M: \psi_s^{-1}(M)
\vert}\sum_{\{mM\in N/M\,:\,\psi_s^{-1}(m)M=nM\}}Y_{mM}\ \text{ for
$n\in N$}.
\label{covariance_with_Y}
\end{equation}
Saying that $U$ is a minimal dilation of $V$ says that $H\subset\H$ is
invariant for $U_s$, that $U_s|_H=V_s$, and that $\bigcup_{s\in S}
U_s^*H$ is dense in $\H$. To see that there is a well-defined
map $W_n$ satisfying (\ref{defW}), we have to prove that
if $U_s^*h=U_t^*k$ for some $s,t\in S$ and $h,k\in H$, then
\begin{equation}\label{welldefW}
U_s^*Y_{\psi_s(n)M}h=U_t^*Y_{\psi_t(n)M}k\ \text{ for every $n\in N$.}
\end{equation}

In proving (\ref{welldefW}), the idea is to move out to some larger
$U_r^*H$ where we can compare the two sides. Formally, we note that
because $S$ is Ore, there are $u,v\in S$ such that $us=vt$, and then
\begin{equation}\label{Vequal}
V_uh=U_{u}h=U_{us}(U_s^*h)=U_{vt}(U_t^*k)=U_vk=V_vk.
\end{equation}
Next, we fix $n\in N$, and note that
\[
\{mM\in N/M:\psi_u^{-1}(m)M=nM\}=\{\psi_u(n)pM: p\in N\text{ and
}\psi_u^{-1}(p)M=M\}.
\]
Thus the covariance relation (\ref{covariance_with_Y}) gives
\begin{align*}
U_s^*(Y_{\psi_s(n)M}h)
&=U_{us}^*U_u(Y_{\psi_s(n)M}h)=U_{us}^*V_u(Y_{\psi_s(n)M}V_u^*V_uh)\\
&=U_{us}^*\Big(\frac{1}{\vert M: \psi_u^{-1}(M)
\vert}\sum_{\{mM\,:\,\psi_u^{-1}(m)M=\psi_s(n)M\}}Y_{mM}\Big)V_uh\\
&=U_{us}^*Y_{\psi_{us}(n)M}\Big(\frac{1}{\vert M: \psi_u^{-1}(M)
\vert}\sum_{\{pM\,:\,\psi_u^{-1}(p)M=M\}}Y_{pM}\Big)V_uh\\
&=U_{us}^*Y_{\psi_{us}(n)M}(V_uY_MV_u^*)V_uh
=U_{us}^*Y_{\psi_{us}(n)M}(V_uV_u^*)V_uh\\
&=U_{us}^*Y_{\psi_{us}(n)M}(V_uh).
\end{align*}
A similar computation shows that
$U_t^*(Y_{\psi_t(n)M}k)=U_{vt}^*Y_{\psi_{vt}(n)M}(V_vk)$, which by
(\ref{Vequal}) is also equal to $U_{us}^*Y_{\psi_{us}(n)M}(V_uh)$. This
gives (\ref{welldefW}).

We have now shown that there is a well-defined map $W_n$ on $\bigcup
U_s^*H$ satisfying (\ref{defW}). It is easy to check that $W_n$ is linear
and isometric on each $U_s^*H$, and hence it extends to an isometry on
$\H$; since $W_mW_n=W_{mn}$, $W$ is a unitary representation of $N$. For
$h\in H$ and $s,t\in S$, we have
\[
U_sW_nU_s^*(U_t^*h)
=U_s U_{ts}^*Y_{\psi_{ts}(n)M}h
=U_t^*Y_{\psi_{t}(\psi_s(n))M}h
= W_{\psi_s(n)}(U_t^*h),
\]
so that $U_sW_nU_s^*=W_{\psi_s(n)}$. Since every $g\in G$ has the form
$g=s^{-1}t$ for $s,t\in S$, it follows that $U_gW_nU_g^*=W_{\psi_g(n)}$
for every $g\in G$, and hence the pair $(W,U)$ does indeed define a
representation $W\rtimes U$ of $N\rtimes_\psi G$. For $m\in M$ we have
$W_mh=h$ for $h\in H$; thus $H\subset \H^M$, and the density of
$\bigcup U_s^*H$ in $\H$ implies that $W\rtimes U$ is generated by its
$M$-fixed vectors. That $W_nh=Y_{nM}h$ for $h\in H$ is clear from
(\ref{defW}).

We have now proved the first part of Theorem~\ref{Hallviadilation}. For
the converse, we start with a representation $X\rtimes U$ of
$N\rtimes_\psi G$ on $\H$ which is generated by $\H^M$. For $m\in M$ and
$h\in \H^M$, we have $X_m(X_nh)=X_n(X_{n^{-1}mn}h)=X_nh$, so $X_nh\in
\H^M$; thus $X|_{\H^M}$ is a unitary representation of $N$ which
defines a representation $Y$ of $N/M$, as claimed. For $m\in M$ and $s\in
S$, we have $X_m(U_sh)=U_s(X_{\psi_s^{-1}(m)}h)=U_sh$, so $U_s$ also
leaves $\H^M$ invariant, and $V_s:=U_s|_{\H^M}$ defines an isometric
representation of $S$. Since $\H$ is generated by $\H^M$, the space
$\H_0:=\bigcup_{(n,g)\in N\rtimes\,G} U_gX_n\H^M$ is dense in $\H$;
since both $X$ and $U|_S$ leave $\H^M$ invariant and $G=S^{-1}S$, we
have $\H_0=\bigcup_{s\in S}U_s^*\H^M$. Thus $U$ is a minimal unitary
dilation of $V$.

To verify covariance of $(\pi_Y,V)$, we need the following standard lemma
about representations of finite groups.

\begin{lem}
Let $R:G_0\to U(H)$ be a representation of a finite group $G_0$ on
a Hilbert space $H$. Then $|G_0|^{-1}\sum_{g\in G_0}R_g$ is
the orthogonal projection on the subspace of $G_0$-fixed vectors.
\label{rep_finite_group}
\end{lem}

\begin{cor}\label{projection_on_H}
Let $s\in S$ and
suppose $F\subset M$ contains one representative of each coset in
$\{mM\in N/M:\psi_s^{-1}(m)M=M\}$.  Then $|M:\psi_s^{-1}(M)|^{-1}
\sum_{m\in F}X_{m}|_{U_s^*\H^M}$ is the
projection $P_s$ of $U_s^*\H^M$ on $\H^M$.
\end{cor}

\begin{proof}
For $m\in M$, $s\in S$ and $h\in \H^M$, we have
$X_m(U_s^*h)=U_s^*(X_{\psi_s(m)}h)$, which belongs to $U_s^*\H^M$ because
$\H^M$ is invariant for $X$. Similarly, we have
\[
X_{\psi_s^{-1}(m)}U_s^*h
=U_s^*X_{\psi_s(\psi_s^{-1}(m))}h =U_s^*X_mh
=U_s^*h,
\]
so
$X|_{\psi_s^{-1} (M)} \equiv 1$ on  $U_s^*\H^M$.
Thus $X$ induces a representation $R$ of $M/\psi_s^{-1}(M)$ on
$U_s^*\H^M$. The space of $R$-fixed vectors is $\H^M$, and
hence the result follows from Lemma \ref{rep_finite_group}.
\end{proof}

We verify covariance of $(\pi_Y,V)$ by checking
(\ref{covariance_with_Y}). With the notation of
Corollary~\ref{projection_on_H}, the solutions of
$\psi_s^{-1}(m)M=nM$ in $\psi_s(M)/M$
are $\{\psi_s(nm)M:m\in F\}$.
Thus applying the right-hand side of (\ref{covariance_with_Y}) to $h\in
\H^M$ gives
\begin{align}\label{checkcov42}
\Big(\frac{1}{|M: \psi_{s^{-1}}(M)|}\sum_{m\in F}Y_{\psi_s(nm)M}\Big)h
&=\frac 1{|M :\psi_s^{-1}(M)|} \sum_{m\in F} X_{\psi_s(nm)}h\\
&=\frac 1{|M :\psi_s^{-1}(M)|} \sum_{m\in F}
U_sX_{nm}U_s^*h\notag\\
&=U_sX_n\Big(\frac{1}{|M :\psi_s^{-1}(M)|}\sum_{m\in F} X_m\Big)U_s^*h
\notag\\
&=U_sX_nP_sU_s^*h\notag,
\end{align}
where we used Corollary \ref{projection_on_H} at the last step. Now note,
first, that $P_sU_s^*h\in \H^M$, so $X_nP_sU_s^*h=Y_{nM}P_sU_s^*h$, and,
second, that $P_sU_s^*$ is the adjoint of $V_s=U_s|_{\H^M}$ in $B(\H^M)$.
Thus $U_sX_nP_sU_s^*h=V_sY_{nM}V_s^*h$, and it follows from
(\ref{checkcov42}) that $(\pi_Y,V)$ is a covariant representation of
$(C^*(N/M),S,\alpha)$.

Finally, that we recover $W=X$ from the construction of the first part
follows from Equation (\ref{defW}), the covariance of $(\pi_X,U)$, and the
definition of $Y_{nM}$ as the restriction of $X_n$: for $h\in \H^M$,
\[
W_n(U^*_sh)=U_s^*\pi_Y(\delta_{\psi_s(n)M})h=U_s^*Y_{\psi_s(n)M}h
=U_s^*X_{\psi_s(n)}h=X_n(U_s^*h).
\]

This completes the proof of Theorem~\ref{Hallviadilation}.

\section{The minimal automorphic dilation}\label{secminauto}

Theorem~\ref{Hallviadilation} tells us that we are interested in
dilations of representations of the semigroup dynamical system
$(C^*(N/M),S,\alpha)$. Laca's analysis in \cite[\S2]{L} says that we
should therefore be looking at the \emph{minimal automorphic dilation}
$(B,G,\beta)$ of $(C^*(N/M),S,\alpha)$, which is an automorphic dynamical
system such that $B\rtimes_\beta G$ contains $C^*(N/M)\rtimes_\alpha S$
as a full corner. The minimal dilation $(B,G,\beta)$ is characterised up
to isomorphism by the existence of an embedding $i:C^*(N/M)\to B$ such
that $\beta_s\circ i=i\circ \alpha_s$ for $s\in S$ and $\bigcup_{s\in
S}\beta_s^{-1}\circ i(C^*(N/M))$ is dense in $B$ \cite[Theorem~2.1]{L}.

To identify the minimal automorphic dilation, we build a candidate for
the dilation and apply \cite[Theorem~2.1]{L}. This was previously done
for $(C^*(\Q/\Z),\N^*,\alpha)$ in \cite[\S3]{L}; the
dilation of this system is based on the algebra $C_0(\A_f)$ of
functions on the abelian group
$\A_f$ of finite adeles, and the argument justifying this in \cite{L} uses
the Fourier transform and standard properties of the adeles. Here  our main
task is to show that in the situation of \S\ref{Hallsthm} there is a good
analogue of
$\A_f$; this construction may be of interest in its own right (see \S\ref{exs}).
Because $N/M$ need not be abelian, we avoid the Fourier transform, and work
directly with group $C^*$-algebras.

Because it is Ore, the generating semigroup $S$ directs $G$: $g\leq_r h
\Longleftrightarrow h\in Sg$. (The subscript
$r$ in $\leq_r$ reminds us that the relation $\leq_r$ is
right-invariant.) If $s,t\in S$
satisfy $s\leq_r t$, then $t=rs$ for some $r\in S$, and
$\psi_t^{-1}(M)=\psi_s^{-1}(\psi_r^{-1}(M))\subset \psi_s^{-1}(M)$. Thus
the quotient maps
\[
\pi_t^s:N/\psi_t^{-1}(M) \to N/\psi_s^{-1}(M)
\]
form an inverse system over the directed set $(S,\leq_r)$. We denote
by $\pi^s$ the canonical homomorphism of the inverse limit
$N_\infty:=\varprojlim N/\psi_s^{-1}(M)$
into $N/\psi_s^{-1}(M)$, so that $\pi^s=\pi_t^s\circ\pi^t$ for $s\leq_r
t$. We often
view
$N_\infty$ as the subgroup of the direct product
$\prod_{s\in S}N/\psi_s^{-1}(M)$ consisting of the families
$(x_s)$ satisfying the compatibility relations $\pi_t^s(x_t)=x_s$ for
$s\leq_r t$; this viewpoint shows that, because each $\pi_t^s$ is
surjective, so is each
$\pi^s$. For each
$n\in N$, we can define a compatible family by
$n_s:=n\psi_s^{-1}(M)$, and this gives a canonical homomorphism $j$ of
$N$ into $N_\infty$; the property (\ref{j_injective}) implies that $j$ is
injective.

Since $M$ is normal in $N$, the quotients $M/\psi_s^{-1}(M)$ are
normal in $N/\psi_s^{-1}(M)$, and the inverse limit $K:=\varprojlim
M/\psi_s^{-1}(M)$ is a normal subgroup of $N_\infty$. Because each
$M/\psi_s^{-1}(M)$ is finite, $K$ is a compact topological group. The
open neighbourhoods of the identity in $K$ are the neighbourhood basis
for a unique locally compact group topology on $N_\infty$ (this follows
from
\cite[Theorem~4.5]{HR}, for example). In this topology, $K$ is a compact
open normal subgroup of $N_\infty$, and this has some immediate
consequences:

\begin{lem}\label{cont_on_subgroup}
Suppose a locally compact group $G$ has a compact open normal subgroup
$K$. Then
$G$ is unimodular, and a
homomorphism of $G$ into a topological group is continuous if and
only if it is continuous on $K$.
\end{lem}

\begin{proof}
To define a Haar integral on $G$, we take a Haar measure $\mu$ on $K$, a
family $\{g_i\}$ of coset representatives for $G/K$, and define
\[
\int_G f(g)\,dg=\sum_i \int_K f(g_ik)\,d\mu(k);
\]
it is easy to check that this does not depend on the choice of coset
representatives. If $t\in G$, then $\{tg_i\}$ is another set of coset
representatives, so the integral is left-invariant; since the
automorphism
$k\mapsto t^{-1}kt$ preserves Haar measure on the compact group $K$, and
since $\{g_it\}$ is another set of coset representatives, it is also
right-invariant. To verify the comment on continuity of homomorphisms,
consider a convergent net
$g_i\to g$, and notice that $g_ig^{-1}$ converges to $e$ and is
eventually in
$K$.
\end{proof}

In our situation, it follows that the homomorphism $\pi^e$ is
continuous, and hence induces an isomorphism of the locally compact
quotient
$N_\infty/K$ onto the discrete group $N/M$. Since the cylinder sets
$\{k:\pi^s(k)=m\psi_s^{-1}(M)\}$
form a basis for the topology on $K$, the image $j(M)$ of $M$ under the
embedding $j$ is dense in $K$; since  $\pi^e(j(n))=nM$ and
$\ker \pi^e=K$, every element of $N_\infty$
has the form
$j(n)k$ for some $k\in K$, and it follows easily that
$j(N)$ is dense in $N_\infty$.

\begin{examples}\label{examples_of_N_infty} 
(a) Fix a prime $p$, and take $S=\N$ (under addition),
$G=\Z$, $N=\Z[p^{-1}]$, $\psi_n(r)=p^{-n}r$, and $M=\Z$. Then
$K=\varprojlim \Z/p^l\Z$ is the additive group in the ring $\Z_p$ of
$p$-adic integers, and $N_\infty$ is the additive group in the field
$\Q_p$ of $p$-adic rationals.

(b) We take $S=\N^*$, $G=\Q_+^*$, $N=\Q$,
$\psi_g(r)=g^{-1}r$, and
$M=\Z$, so that the semigroup system is the one which gives the
Bost-Connes Hecke algebra \cite{LR3}.  Since the operation in $\N^*$ is
multiplication, the direction on $\N^*$ is given by $n\leq_r
m\Longleftrightarrow n\,|\,m$. We claim that
$K=:\varprojlim \Z/n\Z$ is the additive group of integer adeles $\mathcal{Z}$,
which is by definition the product $\prod_p
\Z_p$ over all primes $p$, and
$\Q_\infty:=\varprojlim\Q/n\Z$ is the additive group in the ring $\A_f$ of
finite adeles, which is by definition the restricted direct product
\[
\A_f:=\textstyle{\prod_p (\Q_p,\Z_p):=\big\{(x_p)\in\prod_p\Q_p:
x_p\in\Z_p\text{ for all but finitely many $p$}\big\}}
\]
(see \cite[Chapter~5]{RV}, for example).

To check the first claim, note that for each fixed $p$, the map
$(x_n)\mapsto (x_{p^l})$ is a continuous homomorphism of $K=\varprojlim
\Z/n\Z$ into $\Z_p=\varprojlim \Z/p^l\Z$. Thus the map
$\phi:(x_n)\mapsto ((x_{p^l}))$ is a continuous homomorphism of $K$ into
$\mathcal{Z}$. If $\phi((x_n))=0$, and $n\in \N^*$ has prime
factorisation $n=\prod p^{l_p}$, then $x_n\equiv x_{p^{l_p}}\equiv 0\pmod
{p^{l_p}}$ for all $p$, and $x_n\equiv 0\pmod n$; thus $\phi$ is
injective. It is also surjective: if $((x_{p^l}))\in\mathcal{Z}$, and for
$n=\prod p^{l_p}$ we take $x_n$ to be the unique solution (mod $n$) of the
congruences $x_n\equiv x_{p^{l_p}}\pmod {p^{l_p}}$ guaranteed by the
Chinese Remainder Theorem, then $(x_n)$ is a compatible family with
$\phi((x_n))=((x_{p^l}))$. Since the range of $\phi$ is compact,
the inverse is continuous too, and $\phi$ is the required
isomorphism. In fact, because each coordinate map is a ring homomorphism,
$\phi$ is an isomorphism of compact topological rings, and it
carries the image
$j(n)$ of each integer $n$ into the element $(j_p(n))$ of the product
determined by the embeddings
$j_p:\Z\to \Z_p$.

Since $\Q=\bigcup_{m\in \N}m^{-1}\Z$, the group $\Q_\infty$ is the union
of the subgroups
$m^{-1}K$, each of which is isomorphic to $K$ via the maps $x\mapsto mx$.
Similarly,
$\prod_p (\Q_p,\Z_p)$ is the union of the subgroups $m^{-1}\mathcal{Z}$,
where $m^{-1}$ is the inverse of the element $(j_p(m))$ of $\prod\Q_p$;
since $j_p(m)$ is a unit in $\Z_p$ for all but finitely many $p$, the
inverse makes sense in the restricted product. Now we define
$\mu_m:m^{-1}K\to m^{-1}\mathcal{Z}$ by $\mu_m(x):=m^{-1}\phi(mx)$. If
$n$ divides $m$, say $m=kn$, and $x\in n^{-1}K\subset m^{-1}K$, then
$nx\in K$ and the multiplicativity of $\phi$ gives
\begin{align*}
\mu_{kn}(x)&=(kn)^{-1}\phi(knx)=(kn)^{-1}\phi(j(k))\phi(nx)\\
&=(kn)^{-1}(j_p(k))\phi(nx)=n^{-1}\phi(nx)=\mu_n(x).
\end{align*}
So the $\mu_m$ combine to give a well-defined isomorphism $\mu$ of
$\Q_\infty$ onto $\A_f=\bigcup_m m^{-1}\mathcal{Z}$. Because
$\mu|_K=\phi$ is a homeomorphism, so is $\mu$, and we have proved the
second claim.
\end{examples}

The automorphic dilation which we seek will be based on the
$C^*$-algebra of the locally compact group $N_\infty$. To describe the
action, recall that the automorphisms $\psi_t$  are isomorphisms
of
$\psi_{st}^{-1}(M)$ onto $\psi_{s}^{-1}(M)$ for every $s\in S$, and
hence induce isomorphisms of $N/\psi_{st}^{-1}(M)$ onto
$N/\psi_{s}^{-1}(M)$. To avoid excessive twiddling, we use $\psi_t$
to denote these and other homomorphisms induced by $\psi_t$.

\begin{lem}\label{jexists}
There is an action $\theta$ of $G$ by continuous automorphisms
of $N_\infty$ such that
$\pi^s\circ \theta_t=\psi_t\circ \pi^{st}$ for every $s,t\in S$, and we
then have $j\circ \psi_t=\theta_t\circ j:N\to N_\infty$.
\end{lem}

\begin{proof}
For fixed $t$, we apply the universal property of the inverse limit to
the family
$\psi_t\circ \pi^{st}:N_\infty\to N/\psi_s^{-1}(M)$ to obtain an
endomorphism
$\theta_t$ of $N_\infty$ such that $\pi^s\circ \theta_t=\psi_t\circ
\pi^{st}$ for every $s$.  Applying
the same universal property to
$\psi_t^{-1}\circ \pi^s$ gives an inverse for $\theta_t$, so it is an
automorphism. Since the homomorphisms $\psi_t^{-1}\circ \pi^s$ are
continuous on $K$ and take values in $M/\psi_{st}^{-1}(M)$, the universal
property of the topological inverse limit $K=\varprojlim
M/\psi_s^{-1}(M)$ implies that $\theta_t^{-1}$ is continuous on $K$;
since $K$ is compact, this implies that  $\theta_t^{-1}$ is a
homeomorphism of $K$ onto its image $(\pi^t)^{-1}(e)$. Because both $K$
and $(\pi^t)^{-1}(e)$ are compact open subgroups of $N_\infty$, this is
enough to ensure that $\theta_t^{-1}$ and $\theta_t$ are continuous. To
see that $j\circ \psi_t=\theta_t\circ j$, check that $\pi^s(j\circ
\psi_t(n))=\pi^s(\theta_t\circ j(n))$ for all $s\in S$ and $n\in
N$.
\end{proof}

The action $\theta$ of $G$ on $N_\infty$ induces an action
$\theta_*:G\to\Aut
C^*(N_\infty)$, which will be the action in our
dilated system. We now describe the embedding of
$C^*(N/M)$ in $C^*(N_\infty)$:

\begin{lem}\label{defi}
Choose the Haar measure $\mu$ on $N_\infty$ such that $\mu(K)=1$.
Then $\chi_K\in L^1(N_\infty, d\mu)\subset C^*(N_\infty)$ is a projection
in $C^*(N_\infty)$, and there is a unital embedding
$i$ of
$C^*(N/M)$ in
$\chi_KC^*(N_\infty)\chi_K$ such that
$i(\delta_{nM})=\chi_{j(n)K}$ for all $n\in N$.
\end{lem}

\begin{proof}
For $x,y,z\in N_\infty$, we have
\begin{align*}
\chi_{xK}*\chi_{yK}(z)
&=\int_{N_\infty}\chi_{xK}(w)\chi_{yK}(w^{-1}z)\,d\mu(w)\\
&=\int_{N_\infty}\chi_{K}(x^{-1}w)\chi_{yK}(w^{-1}z)\,d\mu(w)\\
&=\int_{K}\chi_{yK}(w^{-1}x^{-1}z)\,d\mu(w).
\end{align*}
Now because $K$ is normal in $N_\infty$, we have
\begin{align*}
w^{-1}x^{-1}z\in yK\text{ and }w\in K
&\Longleftrightarrow Kx^{-1}z=yK\text{ and }w\in K\\
&\Longleftrightarrow x^{-1}zK=yK\text{ and }w\in K\\
&\Longleftrightarrow zK=xyK\text{ and }w\in K.
\end{align*}
Thus
\begin{align*}
\chi_{xK}*\chi_{yK}(z)
&=\begin{cases}
\int_K \chi_K(z)\,d\mu(z)&\text{if $zK=xyK$}\\
0&\text{otherwise}
\end{cases}\\
&=\chi_{xyK}(z).
\end{align*}
This implies in particular that $\chi_K^2=\chi_K$ in $C^*(N_\infty)$.
Since $N_\infty$ is unimodular, we trivially have $\chi_K^*=\chi_K$, so
$\chi_K$ is a projection in $C^*(N_\infty)$. The above calculation shows
that the elements $\chi_{xK}$ lie in the corner
$\chi_KC^*(N_\infty)\chi_K$, and that the map $xK\mapsto \chi_{xK}$ is a
homomorphism of $N_\infty/K$ into the unitary group of the corner.
Composing with the isomorphism $nM\mapsto j(n)K$ of $N/M$ onto
$N_\infty/K$ gives a homomorphism on $N/M$, and the integrated form $i$
of this homomorphism has the required property.
\end{proof}

\begin{thm}\label{iddilation}
The system $(C^*(N_\infty),G,\theta_*)$ is a minimal automorphic dilation
of the semigroup system $(C^*(N/M),S,\alpha)$: the
homomorphism $i:C^*(N/M)\to C^*(N_\infty)$ of Lemma~\ref{defi} is
injective and satisfies
$i\circ \alpha_s=(\theta_*)_s\circ i$ for $s\in S$, and the union of the
subalgebras
$(\theta_*)_s^{-1}(i(C^*(N/M)))$ is dense in $C^*(N_\infty)$.
\end{thm}

We
need to know how $\theta_*$ acts on the dense $*$-subalgebra
$L^1(N_\infty)$:

\begin{lem}\label{thetaonL1}
For $f\in L^1(N_\infty)\subset C^*(N_\infty)$ and $s\in S$, we have
\begin{equation}\label{theta*}
(\theta_*)_s(f)=|\theta_s(K):K|^{-1}(f\circ \theta_s^{-1})
=|M:\psi_s^{-1}(M)|^{-1}(f\circ \theta_s^{-1}).
\end{equation}
\end{lem}

\begin{proof}
If $\mu$ is a Haar measure on $N_\infty$, so is $E\mapsto
\mu(\theta_s(E))$, and hence there exists $c\in (0,\infty)$ such that
$\mu\circ\theta_s=c\mu$; plugging in $E=K$ shows that
$c=|\theta_s(K):K|$. Now $\theta_s(K)=\varprojlim
\psi_s(M)/\psi_t^{-1}(M)$ and $K=\varprojlim M/\psi_t^{-1}(M)$, and each
$M/\psi_t^{-1}(M)$ has the same index in $\psi_s(M)/\psi_t^{-1}(M)$, so
we have
\[
|\theta_s(K):K|=|\psi_s(M):M|=|M:\psi_s^{-1}(M)|
\]
(this follows from
\cite[Lemma~3.6]{LPR}, for example). This implies that
\begin{equation}\label{inttheta}
\int f\circ\theta_s \,d\mu=|\theta_s(K):K|^{-1}\int f\,d\mu
=|M:\psi_s^{-1}(M)|^{-1}\int f\,d\mu
\end{equation}
for $f=\chi_E$, and we can extend (\ref{inttheta})  to
$f\in L^1(N_\infty)$ by the usual bootstrap
arguments. Applying (\ref{inttheta}) to $(\theta_*)_s(f)=\int
f(x)\delta_{\theta_s(x)}\,d\mu(x)\in C^*(N_\infty)$ gives (\ref{theta*}).
\end{proof}

\begin{proof}[Proof of Theorem~\ref{iddilation}]
To see that $i$ is injective, let $\pi=\pi_W$ be a faithful
representation of $C^*(N/M)$. Then for $n\in N$, we have
\begin{align*}
\pi_{W\circ\pi^e}(i(\delta_{nM}))
&=\int_{N_\infty} \chi_{j(n)K}(x)W\circ\pi^e(x)\,d\mu(x)\\
&=\int_K W\circ\pi^e(j(n)x)\,d\mu(x)\\
&=W_{nM}\Big(\int_K W\circ\pi^e(x)\,d\mu(x)\Big),
\end{align*}
which is just $W_{nM}$ because $K=\ker \pi^e$ and $\mu(K)=1$. Thus the
representation $\pi_W$ factors as $\pi_{W\circ\pi^e}\circ i$, and $i$ must
be injective.

To verify that $i\circ \alpha_s=(\theta_*)_s\circ i$, let $nM\in N/M$.
Then
\[
i\circ\alpha_s(\delta_{nM})
=|M:\psi_s^{-1}(M)|^{-1}\sum_{\{mM\,:\,\psi_s^{-1}(mM)=nM\}}\chi_{j(m)K}.
\]
The cosets $mM$ such that $\psi_s^{-1}(mM)=nM$ are disjoint, with union
$\psi_s(nM)\subset N$; the corresponding cosets $j(m)K=(\pi^e)^{-1}(mM)$
in $N_\infty$ are also disjoint, with union $\theta_s(j(n)K)$. Thus from
Lemma~\ref{thetaonL1} we have
\begin{align*}
i\circ\alpha_s(\delta_{nM})
&=|M:\psi_s^{-1}(M)|^{-1}\chi_{\theta_s(j(n)K)}
=|M:\psi_s^{-1}(M)|^{-1}\chi_{j(n)K}\circ \theta_s^{-1}\\
&=(\theta_*)_s(\chi_{j(n)K})
=(\theta_*)_s\circ i(\delta_{nM}),
\end{align*}
which implies that $i\circ \alpha_s=(\theta_*)_s\circ i$ on all of
$C^*(N/M)$.

For the minimality, we use Lemma~\ref{thetaonL1} again to see that
\begin{align}\label{minimality}
(\theta_*)_s^{-1}(i(\delta_{nM}))&=(\theta_s^{-1})_*(\chi_{j(n)K})\\
&=\vert K:\theta_s^{-1}(K) \vert \chi_{j(n)K}\circ \theta_s\notag\\
&=\vert K:\theta_s^{-1}(K) \vert
\chi_{(\pi^s)^{-1}(\psi_s^{-1}(nM))}.\notag
\end{align}
Since the cylinder sets $\{(\pi^s)^{-1}(z):z\in M/\psi_s^{-1}(M)\}$ form
a basis of compact open sets for the topology on $K$, the sets
\[
\{(\pi^s)^{-1}(\psi_s^{-1}(nM)):n\in N\}
\]
form a basis of compact open sets for the topology on $N_\infty$. Thus
the functions (\ref{minimality}) span a dense subspace of
$C_c(N_\infty)\subset L^1(N_\infty)$, and hence also of $C^*(N_\infty)$.
\end{proof}

\section{The main theorem}

We now identify the category of representations of $N\rtimes_\psi G$
which are generated by their $M$-fixed vectors. We retain the notation
of the previous sections. In particular, we recall from
Lemma~\ref{jexists} that the embedding
$j:N\to N_\infty$  intertwines the actions $\psi$
and $\theta$ of $G$, and hence induces an
embedding of $N\rtimes_\psi G$ as a dense subgroup of the locally
compact group $N_\infty\rtimes_\theta G$. We identify  $N\rtimes_\psi G$
with its image in $N_\infty\rtimes_\theta G$.

\begin{thm}\label{idcategory}
Restriction from $N_\infty\rtimes_\theta G$ to $N\rtimes_\psi G$ is an
equivalence between the categories of continuous unitary representations
of $N_\infty\rtimes_\theta G$ and the category of unitary
representations of $N\rtimes_\psi G$ which are generated by their
$M$-fixed vectors. Both categories are equivalent to the category of
representations of 
$C^*(N/M)\rtimes_\alpha S$.
\end{thm}

\begin{proof}
Since $N\rtimes_\psi G$
is dense in $N_\infty\rtimes_\theta G$, an operator intertwines two
representations of
$N_\infty\rtimes_\theta G$ if and only if it intertwines their
restrictions to $N\rtimes_\psi G$. So for the first part, it suffices to
prove that every restriction is generated by its $M$-fixed vectors, and
that every representation of $N\rtimes_\psi G$ which is generated by its
$M$-fixed vectors extends to a representation of $N_\infty\rtimes_\theta
G$.

Suppose $W\rtimes U$ is a continuous unitary representation of
$N_\infty\rtimes_\theta G$ on $\H$. Since
$C^*(N_\infty)=\overline{\bigcup_{s\in
S}(\theta_*)_s^{-1}(i(C^*(N/M)))}$ acts nondegenerately, the
representation is generated by the vectors in $\pi_W(i(1_{C^*(N/M)}))\H$.
But
$i(1_{C^*(N/M)})=\chi_K$, so for $m\in M$, we have
\begin{align*}
W_{j(m)}\pi_W(i(1_{C^*(N/M)}))&=W_{j(m)}\Big(\int_K W_k\,d\mu(k)\Big)\\
&=\int W_{j(m)k}\,d\mu(k)\\ 
&=\pi_W(i(1_{C^*(N/M)})),
\end{align*}
and hence $W\rtimes U|_M=W\circ j|_M$ fixes $\pi_W(i(1_{C^*(N/M)}))\H$.

Next suppose that $X\rtimes U$ is a unitary representation of
$N\rtimes_\psi G$ which is generated by its space $\H^M$ of
$M$-fixed vectors. From Theorem~\ref{Hallviadilation}, we
know that
$U$ is a minimal dilation of $V:=U|_{\H^M}$, and if we set
$Y_{nM}:=X_n|_{\H^M}$, then
$(\pi_Y,V)$ is a covariant representation of $(C^*(N/M),S,\alpha)$ on
$\H^M$. Because $(C^*(N_\infty),G,\theta_*)$ is a minimal dilation
of $(C^*(N/M),S,\alpha)$, we know from \cite[Lemma~2.3]{L} that there is
a unique representation $\pi$ of $C^*(N_\infty)$ on $\H$ such that
$(\pi, U)$ is a covariant representation of $(C^*(N_\infty),G,\theta_*)$
and $\pi\circ i|_{\H^M}=\pi_Y$. Let $W$ be the continuous unitary
representation of $N_\infty$ such that $\pi=\pi_W$; we claim that
$W\rtimes U$ is the required continuous unitary representation of
$N_\infty\rtimes_\theta G$. It suffices to see that $W|_N=X$, or
equivalently that $W\circ j=X$. Let $n\in N$, and note that
\[
\pi\circ i(\delta_{nM})=\pi(\chi_{j(n)K})=W_{j(n)}\pi(\chi_K)
=W_{j(n)}\pi\circ i(1_{C^*(N/M)}).
\]
The relation $\pi\circ i|_{\H^M}=\pi_Y$ implies that $\pi\circ
i(1_{C^*(N/M)})$ is the identity on $\H^M$, and hence that
\[
W_{j(n)}|_{\H^M}=\pi\circ i(\delta_{nM})|_{\H^M}=Y_{nM}.
\]
Thus for $s\in S$ and $h\in \H^M$, we have
\begin{align*}
W_{j(n)}U^*_sh&=U_s^*W_{\theta_s(j(n))}h=U_s^*W_{j(\psi_s(n))}h
=U_s^*Y_{\psi_s(n)M}h\\
&=U_s^*X_{\psi_s(n)}h=X_nU_s^*h,
\end{align*}
which implies $W_{j(n)}=X_n$ because $\bigcup U_s^*\H^M$ is dense in $\H$.

The last assertion follows from Theorem~\ref{iddilation} and
Theorem~\ref{functorialLaca}.
\end{proof}

\section{Examples}\label{exs}

As we mentioned in the introduction, our construction of $N_\infty\rtimes G$
as an inverse limit seems to be a little different from that used by Schlichting
and other authors, and may have some advantages. The presentation of the
additive adeles
$\A_f$ as the inverse limit 
$\varprojlim\Q/n\Z$ over the directed set $\N^*$, for example, is a little
unusual. Here we shall see how it quickly gives the  self-duality of $\A_f$,
and use similar arguments to produce a family of self-dual locally compact
abelian groups from the examples considered in \cite{B} and \cite{LarR}.

\begin{example}\label{self_duality_of_adeles}
We use the descriptions $\A_f=\varprojlim\Q/n\Z$ and
$\mathcal{Z}=\varprojlim \Z/n\Z$ from Example 
\ref{examples_of_N_infty}. We begin by recalling that the pairings $\langle
k,r\rangle:=\exp(2\pi ikr)$ of $\Z/n\Z$ with ${\frac 1n \Z}/\Z$ induce
isomorphisms of $\Z/n\Z$ onto $({\frac 1n \Z}/\Z)^\wedge$. These pairings are
compatible with the bonding maps $\Z/m\Z\to\Z/n\Z$, which are given by
reduction mod $n$ when $n$ divides $m$, and with the inclusions of ${\frac 1n
\Z}/\Z$ in ${\frac 1m \Z}/\Z$; thus they induce a continuous isomorphism $\Phi$
of
$\mathcal{Z}$ onto $(\Q/\Z)^\wedge=(\varinjlim {\frac
1n\Z}/\Z)^\wedge=\varprojlim ({\frac 1n \Z}/\Z)^\wedge$ such that
\begin{equation}\label{isoonZ}
\Phi(x)(r)=\langle \pi^n(x),r\rangle = e^{2\pi i\pi^n(x)r}\ \text{ for $r\in
\textstyle{{\frac 1n \Z}}/\Z\subset \Q/\Z$}.
\end{equation}

Now for $x=(x_n)_{n\in \N}\in \A_f$, let $d_x$ be the denominator of
$\pi^0(x)$ in $\Q/\Z$. A straightforward calculation shows that 
\begin{equation}
\langle x, y \rangle:=e^{2\pi i \pi^n(x)\pi^n(y)}
\label{pairing_on_adeles}
\end{equation}  
is the same for each $n$ which is divisible by both $d_x$ and $d_y$. For fixed
$x\in\A_f$, the map $\langle x, \cdot \rangle:\A_f\to \T$ is a
homomorphism. Since $\langle x, \cdot \rangle$ is identically $1$ on 
$d_x\mathcal{Z}$, it is continuous on this compact open subgroup of
$\A_f$, and hence is continuous on $\A_f$ by Lemma \ref{cont_on_subgroup}.  
It follows that $\Psi:x\mapsto
\langle x,\cdot  \rangle$ is a homomorphism of $\A_f$ into 
$\widehat\A_f$. Note that $\Psi$ takes $\mathcal{Z}$ onto
$\mathcal{Z}^\perp\cong(\A_f/\mathcal{Z})^\wedge\cong (\Q/\Z)^\wedge$, and
coincides on $\mathcal{Z}$ with the continuous isomorphism $\Phi$ of
(\ref{isoonZ}). The induced homomorphism $q$ of $\Q/\Z\cong \A_f/\mathcal{Z}$
into
$\widehat\A_f/\mathcal{Z}^\perp\cong\widehat{\mathcal{Z}}$ is given by 
$q(r)(z)=\langle z, r\rangle=\exp{(2\pi i \pi^n(z)r)}$ for $r\in \frac 1n 
\Z/\Z$ and $z\in \mathcal{Z}$; it is therefore the dual of the isomorphism in
(\ref{isoonZ}), and in particular is an isomorphism. We deduce, first, that
$\Psi$ is a group isomorphism, and, second, that both $\Psi$ and $\Psi^{-1}$ are
continuous on compact open subgroups, and hence continuous everywhere.
\end{example}

\begin{example}\label{FM_case}
Let $F$ and $M$ be commuting matrices in $GL_d(\Z)$ such that
$\det F\neq 1$, $\det M\neq 1$ and
$(\det F, \det M)=1$. Let $S=\N^2$,
$N=\bigcup_{m, n}F^{-m}M^{-n}\Z^d\subset \Q^d$, and define $\psi:\Z^2 \to
\Aut N$  by $\psi_{k, l}=F^{-k}M^{-l}$. The subgroup $\Z^d$ of $N$ satisfies
$\Z^d\subset \psi_{m,n}(\Z^d)$, Equation~(\ref{j_injective}) and
\[
\vert \Z^d : F^{m}M^{n}\Z^d \vert=(\det F)^m(\det M)^n<\infty \ \text{ for all
$(m,n)\in \N^2$}
\]
(see \cite[\S 4.4]{B} and \cite[Example 3.6]{LarR}). The construction of
\S\ref{secminauto} yields
$$
N_\infty=\varprojlim N/(F^mM^n\Z^d)\  \text{ and }\ K=\varprojlim
\Z^d/(F^mM^n\Z^d),
$$ 
with $N_\infty /K\cong N/\Z^d=\varinjlim\,(F^{-m}M^{-n}\Z^d)/\Z^d.$  Let
$N_\infty^t$ and $K^t$ be the groups obtained by applying the same construction
to the transposes
$F^t$ and $M^t$.

Following the previous example, we fix $(m,n)\in \N^2$, and note that the
usual pairing $(w,k)\mapsto \exp(2\pi i w\cdot k)$ of ${\mathbb R}^d/\Z^d$ with
$\Z^d$ implements an isomorphism of
$\Z^d/(F^mM^n\Z^d)$ onto $((F^t)^{-m}(M^t)^{-n}\Z^d)/\Z^d)^\wedge$; 
these isomorphisms combine to give a continuous  isomorphism $\Phi$ of $K$ onto
$(N^t/\Z^d)^\wedge$. Now for $x\in N_\infty$, we choose the smallest $(m_0,n_0)$
such that
$\pi^0(x)\in (F^{-m_0}M^{-n_0}\Z^d)/\Z^d$. The number
$f(x, y):=\langle\pi^{m, n}(x), \pi^{m, n}(y) \rangle $
is constant for $(m, n)\geq (m_0,n_0)$, and $x\mapsto f(x,\cdot)$
induces a homomorphism $\Psi$ of $N_\infty$  into $(N^t_\infty)^\wedge$ which
restricts to the isomorphism
$\Phi$ of $K$ onto
$(K^t)^\perp\cong(N^t_\infty/K^t)^\wedge\cong (N^t/\Z^d)^\wedge$. The argument
of the previous example shows that 
$\Psi$ is an isomorphism of the topological group $N_\infty$ onto
$(N^t_\infty)^\wedge$. 

The construction of the previous paragraph shows  in particular that the
locally compact abelian group
$N_\infty$ is self-dual whenever $F=F^t$ and
$M=M^t$.
\end{example}

\appendix
\section{Dilation and induction}

Let $\alpha:S\to \End A$ be an action of an Ore semigroup $S$ by
endomorphisms of a $C^*$-algebra $A$. For simplicity, we assume that $A$
has an identity, but not that the $\alpha_s$ preserve it. (If analogues
of these results for endomorphic actions on nonunital algebras are
required, it should be possible to obtain them using the techniques of
\cite[\S4]{PRY}.) We consider the minimal automorphic dilation
$(B,G,\alpha)$ of $(A,S,\alpha)$, and write $i$ for the embedding of $A$
in $B$.

The crossed product $B\rtimes_\beta G$ is generated by a copy
of $B$ and a copy $\{u_s:s\in G\}$ of $G$ in $UM(B\rtimes_\beta G)$; elements of
the form
$\{bu_s: b\in B, s\in G\}$ span a dense subspace of $B\rtimes_\beta
G$. We denote by $p=i(1_A)$ the image of the identity of $A$ in $B$.
Then for $s\in S$,
\[
u_sp=\beta_s(i(1_A))u_s=i(\alpha_s(1_A))u_s= i(1_A)i(\alpha_s(1_A))u_s=pu_sp,
\]
so $v_s:=u_sp$ defines an isometric representation $v$ of $S$ in
$p(B\rtimes_\beta G)p$; the pair $(i, v)$ is a covariant representation
of $(A, S, \alpha)$ in $p(B\rtimes_\beta G)p$. It is shown in
\cite[Theorem~2.4]{L} that $p$ is a full projection in $M(B)$, and that
$i\times v$ is an isomorphism of $A\rtimes_\alpha S$ onto
$p(B\rtimes_\beta G)p$. We identify $A\rtimes_\alpha
S$ with the corner $p(B\rtimes_\beta G)p$ using $i\times v$.

The crucial step in the proof of
\cite[Theorem~2.4]{L} is what Laca calls the \emph{dilation-extension}
of a covariant representation $(\pi,V)$ to a covariant representation
$(\rho,U)$  of $(B,G,\beta)$: indeed, if $U:G\to U(\H)$ is any minimal
unitary dilation of $V$, then there is a compatible representation $\rho$
of $B$ on $\H$
\cite[Lemma~2.3]{L}.  To undo
the extension-dilation process, we restrict $U_s$ to the invariant
subspace
$\rho(p)\mathcal{H}_U$, and compress the representation $\rho \circ i$
to obtain a representation of
$A$ on $\rho(p)\mathcal{H}$; we call the resulting covariant
representation of $(A,S,\alpha)$ the
\emph{restriction-compression} of $(\rho,U)$. Restriction-compression is a
functor
$RC$ from the category of covariant representations of
$(B, G, \beta)$ to the category of covariant representations of $(A, S,
\alpha)$, or, equivalently, from the representations of $B\rtimes_\beta
G$ to the representations of $A\rtimes_\alpha S$.

The main result of this Appendix is the following strengthening of
\cite[Theorem~2.4]{L}.

\begin{thm}\label{functorialLaca} Let $\alpha$ be an action of an Ore
semigroup
$S$ by  injective endomorphisms of a unital $C^*$-algebra $A$, and let
$(B, G,
\beta)$ be a  minimal automorphic dilation of $(A, S,\alpha)$. Then
Laca's dilation-extension construction implements an equivalence between
the category of representations of $A\rtimes_\alpha S$  and the category
of representations of $B\rtimes_\beta G$; the inverse is given by
restriction-compression.
\label{main_appendix}
\end{thm}

To verify that Laca's construction defines a functor, we recognise it as
the inducing construction associated to the imprimitivity bimodule
$X:=(B\rtimes_\beta G)p$, and {\it define} the dilation-extension of a
representation $\pi\times V$ of $A\rtimes_\alpha S$  to be
$X\dashind(\pi\times V)$; the general theory of Rieffel then says that
$X\dashind$ is an equivalence of categories \cite[Theorem~3.29]{tfb}. To
prove the theorem, therefore, we have to verify that $X\dashind$
implements Laca's dilation-extension, and that the inverse functor, which  is
the induction process $\tilde{X}\dashind$ associated to the dual bimodule
$\tilde{X}$ of
\cite[\S3.2]{tfb}, is naturally equivalent to restriction-compression.

\begin{lem}
Suppose that $(\pi, V)$ is a covariant representation of $(A, S, \alpha)$
on $H$. Let $X\dashind(\pi \times V)$ be the induced representation of
$B\rtimes G$ on $X\otimes_{A\rtimes_\alpha S} H$, which satisfies
$X\dashind(\pi \times V) (d)(x\otimes h)= (d\cdot x)\otimes h$  for
$d\in B\rtimes G$ and $x\in X$. Write $X\dashind(\pi\times V)=\sigma
\times  W$. Then $W$ is a minimal unitary dilation of
$V$, and $\sigma \circ i\vert_{i(1)H}=\pi$.
\label{induction_is_dilation}
\end{lem}

\begin{proof}
Since $\langle p,p\rangle_{A\rtimes\,S}=p$ is the identity in
$A\rtimes_\alpha S$, $\phi(h):=p\otimes h$ defines an isometric embedding
$\phi$ of $H$ in $X\otimes_{A\rtimes\,S} H$. The equation $v_s=pu_sp=u_sp$
implies that $\phi(V_sh)=W_s\phi(h)$, so $W$ is a dilation of $V$. To see
that it is minimal, note that the elements of the form $u_s^*i(a)u_t$ span
a dense subspace of $B\rtimes_\beta G$, and hence the elements of the form
$(u_s^*i(a)u_tp)\otimes h$ span a dense subspace of $X\otimes_{A\rtimes\,
S} H$. But
\[
(u_s^*i(a)u_tp)\otimes h=u_s^*i(a)v_t\otimes
h=u_s^*p\otimes\pi(a)V_th=W_s^*(\phi(\pi(a)V_th))
\]
belongs to $W_s^*(\phi(H))$, so it follows that $\bigcup W_s^*(\phi(H))$
is dense in $X\otimes_{A\rtimes\,S} H$. Finally, for $a\in A$ and $h\in
H$, we have
\[
\sigma(i(a))(\phi(h))=(i(a)p)\otimes h=p\otimes \pi(a)h=\phi(\pi(a)h),
\]
which gives the last observation.
\end{proof}

\begin{remark}
At first sight it might seem strange that the Hilbert space of the
dilation is determined by $V$ alone. However, this also happens for
induced representations of dynamical systems: if $\beta:G\to \Aut B$,
$H$ is a subgroup of $G$, and $(\pi,U)$ is a covariant representation
of $(B,H,\beta)$, then the induced representation
$\Ind_{B\rtimes\,H}^{B\rtimes\,G}(\pi\times U)$ acts in the Hilbert space of
$\Ind_H^GU$ (see \cite{T}).
\end{remark}

Next we show that the functor $\tilde{X}\dashind$ is naturally equivalent
to the  restriction-compression functor. It is helpful
to note that $\flat(dp)\mapsto pd^*$ is an
isomorphism of $\tilde{X}:=\{\flat(x):x\in X\}$ onto $p(B\rtimes_\beta
G)$.

\begin{lem}
For each covariant representation $(\rho, U)$ of $(B, G, \beta)$ on
$\mathcal{H}_U$, the map $\Theta_{\rho, U}:pd\otimes h\mapsto
\rho\times U(pd)h$ extends to a unitary operator
$\tilde{X}\otimes_{B\rtimes\,G}\mathcal{H}_U\to \rho(p)\mathcal{H}_U$
which implements an equivalence between $\tilde{X}\dashind(\rho\times U)$
and $RC(\rho, U)$. The iso\-morphisms $\Theta_{\rho, U}$ implement
a natural equivalence between $\tilde{X}\dashind$ and $RC.$
\label{rc_is_inverse}
\end{lem}

\begin{proof}
The map $\Theta_{\rho, U}$ is isometric because the inner product on
$\tilde{X}\otimes_{B\rtimes\,G}\mathcal {H}_U$ is given by
$$
\bigl(pd\otimes h\mid pd\otimes h\bigr)=
\bigl(\rho\times U(\langle pd, pd \rangle_{B\rtimes G})h\mid
h\bigr)_{\mathcal {H}_U} =
\bigl(\rho\times U(d^*pd )h\mid h\bigr)_{\mathcal {H}_U};
$$
it is onto because $\rho\times U$ is nondegenerate. A routine
calculation shows that
$$
\Theta_{\rho, U}\bigl(\tilde{X}\dashind(\rho\times U)(c) (pd\otimes h)
\bigr)=RC(\rho, U)(c)(\Theta_{\rho, U}(pd\otimes h))
$$
for $c\in A\rtimes S$ and $pd\otimes h\in \tilde{X}\otimes_{B\rtimes\,
G}\mathcal {H}_U$.

To check naturality, suppose $T:\H_U\to \H_W$ intertwines $\rho\times U$
and $\sigma\times W$. The morphism $RC(T)$ is just the restriction of $T$,
so
\[
RC(T)\circ\Theta_{\rho, U}(pd\otimes h)
= T(\rho\times U(pd)h)
= \sigma\times W(pd)(Th)
=\Theta_{\sigma, W}(pd\otimes Th),
\]
which is $\Theta_{\sigma,W}\big(\tilde{X}\dashind(T)(pd\otimes h)\big)$.
\end{proof}

Theorem \ref{main_appendix} now follows from Rieffel's general result,
Lemma~\ref{induction_is_dilation}, and Lemma~\ref{rc_is_inverse}.

\end{document}